\documentclass[12pt]{iopart}

\usepackage{graphicx,amscd,amsmath,amssymb,verbatim}
\usepackage[dvips]{hyperref}
\usepackage{textcomp}
\usepackage[OT1,T1]{fontenc}

\begin{document}

\title[Riesz criteria revisited]{The criteria of Riesz, Hardy-Littlewood et al. for the Riemann Hypothesis revisited using similar functions}

\author{Stefano Beltraminelli and Danilo Merlini}

\address{CERFIM, Research Center for Mathematics and Physics, PO Box 1132, 6600 Locarno, Switzerland}
\address{ISSI, Institute for Scientific and Interdisciplinary Studies, 6600 Locarno, Switzerland}

\eads{\mailto{stefano.beltraminelli@ti.ch}, \mailto{merlini@cerfim.ch}}

\begin{abstract}
The original criteria of Riesz and of Hardy-Littlewood concerning the truth of the Riemann Hypothesis (RH) are revisited and further investigated in light of the recent formulations and results of Maslanka and of Baez-Duarte concerning a representation of the Riemann Zeta function. Then we introduce a general set of similar functions with the emergence of Poisson-like distributions and we present some numerical experiments which indicate that the RH may barely be true.
\newline
\newline
\newline
Date: 6 January 2006
\end{abstract}

\submitto{\JPA}

\maketitle

\section{Introduction} 

\noindent It is well known that there are many different criteria
for the truth of the Riemann Hypothesis (RH). Some of these are
not directly related to the important high level computations and
developments concerning the non trivial zeros of the Riemann Zeta
function. In fact, at the beginning of the century M. Riesz, and
later G.H. Hardy and J.E. Littlewood (among other important results
in number theory) found a criterion of ``classical type'' for the
truth of the RH. The above criteria are related to some series involving
values of the Zeta function outside the critical strip, i.e. at
integers arguments of the Zeta function [7, 9], and in a numerical
context, very accurate calculations are needed toward a ``possible
kind of verification`` of the RH.

\noindent In the literature important remarks have been given by
leading mathematicians (see for example, those cited in [3]). We
may think that such criteria may have a limited interest since,
with them, one should work outside the critical strip. It is, in
fact, true that in dealing with the above criteria one needs the
use of arguments of the Zeta function outside the critical strip,
and problems of interchange of summations are present. As an example,
in the above criteria, if one uses the formula established by the
authors, one should give a meaning to an integration over the real
line, which exists only for finite intervals. In order to obtain
finite numerical results which give ``satisfactory'' values to the
functions supposed to be equal to the reciprocal of the Zeta function
outside and inside the critical strip, the integration should be
carried out using a special sequence of upper limit of integration
extending to infinity [6]. 

\noindent But lately, there have been new developments and rigorous
results in connection with this kind of problem: first a ``regularization``
of the representations of the Zeta function (a pioneering work by
Maslanka [8]), followed (in particular) by a new rigorous discrete
formulation with theorems concerning the above criteria (the works
of Baez-Duarte [1, 2, 3]).

\noindent In light of these new approaches, we thought that some
of the above criteria deserve still more study, at least in the
direction of some numerical experiments. Thus, we introduce additional
functions containing two parameters, in order to have additional
confidence in the numerical results of the experiments.

\noindent The content of this work is as follows: in Section 2 we
define a general set of functions with two parameters $\alpha$ and
$\beta$ in the spirit of Riesz and of Hardy-Littlewood and then
obtain the discrete ``representation'' of the Zeta function of our
set by means of the two parameter Pochammer's polynomials with their
coefficients. For the reader the discussion of the conditions are
then given in Appendix A and in Appendix B (they follow strictly
the ingenious method of Baez-Duarte for the Riesz case $\alpha =\beta
=2$).

\noindent In Section 3 we then obtain in some ``limit'', a Poisson
distribution for the coefficients $c_{k}$ of the Pochammer's polynomials;
this is useful in the context of the numerical experiments. These
are presented in Section 4 where many various limiting cases are
treated. In the case of increasing values of the parameter $\beta$,
the experiments indicate that the Poisson distribution becomes more
and more exact and the discrete function $c_{k}$ becomes a constant
which can be evaluated. Finally, we present the results of an experiment
carried out in the critical strip with different values of $\mathfrak{R}(
s) $. As a consequence, we may argue that in the context of the
range of validity of the experiments we present, at large and at
low values of {\itshape k}, the RH appears to be barely true.

\section{The model}

\noindent We now consider a set of functions with two parameters
($\alpha$, $\beta$) to obtain $\frac{1}{\zeta ( s) }$. These are
simply an extension of these two cases: the first (with $\alpha
=\beta =2$) introduced and studied by Riesz [9], the second one
(where $\alpha =1$ and $\beta =2$) by Hardy-Littlewood [7].

\noindent Let $\mu ( n) $ be the M\"obius function of argument $n$,
where
\[
\mu ( n) =\begin{cases}
1,&\mathrm{if} \ n=1\text{}    \\
{\left( -1\right) }^{k},&\mathrm{if} \ n \ \mathrm{is} \ \mathrm{a} \ \mathrm{product} \
\mathrm{of} \ k \ \mathrm{distinct} \ \mathrm{primes}   \\
 0,&\mathrm{if} \ n\  \mathrm{contains} \ \mathrm{a} \ \mathrm{square} \\
\end{cases}
\]
\noindent Let $s_{0}=\rho \mathrm{+}\mathrm{i} \mathit{t}$ and {\itshape
s} a complex variable. For $\mathfrak{R}( s) >\rho =1$ one has $\frac{1}{\zeta
( s) }=\sum \limits_{n=1}^{\infty }\frac{\mu ( n) }{n^{s}}$.

\noindent The two-parameters family of functions is given by:
\begin{equation}
\varphi ( s;\alpha ,\beta ) :=\frac{1}{\Gamma ( -\frac{s-\alpha
}{\beta }) }\int _{0}^{\infty }\sum \limits_{n=1}^{\infty }\frac{\mu
( n) }{n^{\alpha }} e ^{-\frac{x}{n^{\beta }}}x^{-\left( \frac{s-\alpha
}{\beta }+1\right) }dx
\end{equation}

\noindent so that expanding the right-hand side in powers of $x$,
we obtain:
\begin{align*}
\varphi ( s) &=\frac{1}{\Gamma ( -\frac{s-\alpha }{\beta }) }\int
_{0}^{\infty }\sum \limits_{n=1}^{\infty }\frac{\mu ( n) }{n^{\alpha
}}\sum \limits_{k=0}^{\infty }\frac{{\left( -1\right) }^{k}x^{k}}{k!n^{\beta
k}}x^{-\left( \frac{s-\alpha }{\beta }+1\right) }dx
\\%
 &=\frac{1}{\Gamma ( -\frac{s-\alpha }{\beta }) }\int _{0}^{\infty
}\psi ( x) x^{-\left( \frac{s-\alpha }{\beta }+1\right) }dx
\end{align*}

\noindent where
\begin{equation}
\psi ( x) =\sum \limits_{k=0}^{\infty }\frac{{\left( -1\right) }^{k}x^{k}}{k!}\frac{1}{\zeta
( \alpha +\beta k) }
\end{equation}

\noindent The function $\psi ( x) $ was introduced by Riesz (case
$\alpha =\beta =2$) and by Hardy-Littlewood (case $\alpha =1, \beta
=2$).\ \ \ 

\noindent If $\psi ( x) \sim \frac{A}{x^{\frac{\alpha -\rho }{\beta
}-\epsilon }}$ for some $\epsilon$ and for large $x$, then\ \ 
\[
\left| \varphi ( s) \right| \leq \left| \frac{1}{\Gamma ( -\frac{s-\alpha
}{\beta }) }\right| \int _{0}^{\infty }\frac{A}{x^{\frac{\alpha
-\rho }{\beta }+\frac{\mathfrak{R}( s) -\alpha }{\beta }+1-\epsilon
}}dx\leq \left| \frac{1}{\Gamma ( -\frac{s-\alpha }{\beta }) }\right|
\int _{0}^{\infty }\frac{A}{x^{1+\frac{\mathfrak{R}( s) -\rho }{\beta
}-\epsilon }}dx
\]

\noindent would exist and would eventually be given by $\frac{1}{\zeta
( s) }$ with $\zeta ( s) \neq 0$ if we choose $\mathfrak{R}( s)
>\rho +\beta  \epsilon $.

\noindent Let $\rho =\frac{1}{2}$. For $\alpha =\beta =2$ we have:\ \ 
\[
\psi ( x) =\frac{A}{x^{3/4-\epsilon }}
\]

\noindent and for $\alpha =1, \beta =2$: 
\[
\psi ( x) =\frac{A}{x^{1/4-\epsilon }}
\]

\noindent On the other hand expanding (1) in a similar way, we have
that:
\begin{align*}
\varphi ( s) &=\frac{1}{\Gamma ( -\frac{s-\alpha }{\beta }) }\int
_{0}^{\infty }\sum \limits_{n=1}^{\infty }\frac{\mu ( n) }{n^{\alpha
}} e ^{x( 1-\frac{1}{n^{\beta }}) } e ^{-x}x^{-\left( \frac{s-\alpha
}{\beta }+1\right) }dx
\\%
 &=\frac{1}{\Gamma ( -\frac{s-\alpha }{\beta }) }\sum \limits_{k=0}^{\infty
}\sum \limits_{n=1}^{\infty }\frac{\mu ( n) }{n^{\alpha }}{\left(
1-\frac{1}{n^{\beta }}\right) }^{k}\int _{0}^{\infty }\frac{1}{k!}x^{k-\frac{s-\alpha
}{\beta }-1} e ^{-x}dx
\\%
 &=\frac{1}{\Gamma ( -\frac{s-\alpha }{\beta }) }\sum \limits_{k=0}^{\infty
}\sum \limits_{n=1}^{\infty }\frac{\mu ( n) }{n^{\alpha }}{\left(
1-\frac{1}{n^{\beta }}\right) }^{k}\frac{1}{k!}\Gamma ( k-\frac{s-\alpha
}{\beta }) 
\\%
 &=\frac{1}{\Gamma ( -\frac{s-\alpha }{\beta }) }\sum \limits_{k=0}^{\infty
}c_{k}\prod \limits_{r=1}^{k}\left( 1-\frac{\frac{s-\alpha }{\beta
}+1}{r}\right) \Gamma ( -\frac{s-\alpha }{\beta }) 
\end{align*}

\noindent Thus:
\begin{equation}
\varphi ( s) =\sum \limits_{k=0}^{\infty }c_{k}P_{k}( \frac{s-\alpha
}{\beta }+1) 
\end{equation}

\noindent where $P_{k}( x) :=\prod \limits_{r=1}^{k}(1-\frac{x}{r})$
are the Pochhammer polynomials and the functions:
\begin{equation}
c_{k}:=\sum \limits_{n=1}^{\infty }\frac{\mu ( n) }{n^{\alpha }}{\left(
1-\frac{1}{n^{\beta }}\right) }^{k}
\end{equation}

\noindent were already studied by Maslanka [8] and Baez-Duarte [2]
in the special case $\alpha =\beta =2$. Let $\mathfrak{R}( s) >
\rho +\epsilon $ ($\epsilon >0$ and $\rho \in [1,\infty [$). From
a theorem of Baez-Duarte [2], which says that $|P_{k}( s) |\leq
A\cdot k^{-\mathfrak{R}( s) }$ where $A$ is a constant depending
on $|s|$, for large values of $k$ we have that:
\begin{equation}
\left| \varphi ( s) \right| \leq \sum \limits_{k=0}^{\infty }\left|
c_{k}\right| k^{-\left( \frac{\rho +\epsilon -\alpha }{\beta }+1\right)
}
\end{equation}

\noindent In Appendix A we show that if $\alpha >1$ and $\beta
>0$ the following holds unconditionally: 
\[
q_{k}\ll \frac{1}{k^{\frac{\alpha -1}{\beta }}}
\]

\noindent where
\[
q_{k}=\sum \limits_{n=1}^{\infty }\frac{1}{n^{\alpha }}{\left( 1-\frac{1}{n^{\beta
}}\right) }^{k}
\]

\noindent Then we obtain:
\[
\left| \varphi ( s) \right| \leq \sum \limits_{k=0}^{\infty }\frac{1}{k^{\frac{\alpha
-1}{\beta }}}\cdot \frac{A}{k^{\frac{\rho +\epsilon -\alpha }{\beta
}+1}}\leq \sum \limits_{k=0}^{\infty }\frac{A}{k^{\frac{\rho +\epsilon
-1}{\beta }+1}}\leq A\sum \limits_{k=0}^{\infty }\frac{1}{k^{1+\frac{\epsilon
}{\beta }}}<\infty 
\]

\noindent From this it follows that equation (6) below, which gives
${[\zeta ( s) ]}^{-1}$ as:
\begin{equation}
\varphi ( s) =\frac{1}{\zeta ( s) }=\sum \limits_{k=0}^{\infty }c_{k}P_{k}(
\frac{s-\alpha }{\beta }+1) 
\end{equation}

\noindent is valid for $\mathfrak{R}( s) >1$. Now, still from the
theorem of Baez-Duarte [2] i.e. that
\[
\left| P_{k}( \frac{s-\alpha }{\beta }+1) \right| \leq \frac{A}{k^{\frac{\mathfrak{R}(
s) -\alpha }{\beta }+1}}
\]

\noindent and from the hypothesis discussed in Appendix B i.e. supposing
the RH to be true for $\mathfrak{R}( s) >\rho +\epsilon $ ($\epsilon
>0$ and $\rho \in [\frac{1}{2},\infty [$):
\[
\left| c_{k}\right| \ll \frac{B}{k^{\frac{1}{\beta }\left( \alpha
-\rho -\epsilon \right) }}
\]

\noindent then the above series given by (3) converges uniformely.
In fact for $\mathfrak{R}( s) >\rho +\epsilon $ we have: 
\[
\left| \varphi ( s) \right| \leq \sum \limits_{k=0}^{\infty }\frac{B}{k^{\frac{1}{\beta
}\left( \alpha -\rho -\epsilon \right) }}\frac{A}{k^{\frac{\mathfrak{R}(
s) -\alpha }{\beta }+1}}\sim \sum \limits_{k=0}^{\infty }\frac{C}{k^{1+\frac{1}{\beta
}\left( \mathfrak{R}( s) -\rho -\epsilon \right) }}
\]

\noindent Following Baez-Duarte the series $\varphi ( s) $ extends
analytically $\frac{1}{\zeta ( s) }$ to the half plane $\mathfrak{R}(
s) >\rho =\frac{1}{2}$. We have thus obtained for our family of
functions with parameters $\alpha$, $\beta$ that a necessary and
sufficient condition for $\zeta ( s) \neq 0$ in the half plane $\mathfrak{R}(
s) >\rho  $ ($\rho \in [\frac{1}{2},\infty [$) is given by:  
\begin{equation}
\left| c_{k}( \alpha ,\beta ) \right| \leq \frac{\mathrm{const}}{k^{\frac{1}{\beta
}\left( \alpha -\rho -\epsilon \right) }}\ \ \ \ \ \ \ \forall \epsilon
>0, \forall \alpha >1, \forall \beta >0
\end{equation}

\noindent {\itshape Remark 1}

\noindent Instead of using the M\"obius function $\mu$ in $c_{k}$,
one may use (for the numerical computations) the formula involving
values of the Zeta function:
\begin{equation}
c_{k}=\sum \limits_{n=1}^{\infty }\frac{\mu ( n) }{n^{\alpha }}{\left(
1-\frac{1}{n^{\beta }}\right) }^{k}=\sum \limits_{j=0}^{k}{\left(
-1\right) }^{j}\binom{k}{j}\frac{1}{\zeta ( \alpha +\beta j) }\ \ \ \ \ \ 
\end{equation}

\noindent {\itshape Remark 2}

\noindent From the bound above it follows not only theoretically
but also in the context of a numerical analysis that it will be
equally difficult to treat the case $\rho \in [\frac{1}{2},1]$,
for example $\rho =\frac{3}{4}$, as will be the case $\rho =\frac{1}{2}+\epsilon
$ with $\epsilon$ small (see the last experiment concerning the
critical strip where $\rho =\frac{1}{2}$, $\rho =\frac{3}{4}$ and
$\rho =1$. Below in Section 4 we will also treat the case $\alpha
=\frac{7}{2}$. 

\noindent {\itshape Remark 3}

\noindent The condition for the truth of the RH using Riesz and
Hardy-Littlewood functions $\psi ( x) $ is essentially the same
as the one using the discrete function $c_{k}$ with $k\in \mathbb{N}$.
In a previous work [6] independent of the present one (which essentially
uses the Baez-Duarte idea and theorems) some numerical results were
obtained for $\psi (x)$ in the case of the Hardy-Littlewood function
($\alpha =1, \beta =2$) by integration in the $x$-space. The discrete
version using the function $c_{k}$ of the discrete variable $k$
[2, 8] has advantages in the numerical computations which will be
presented below. Before this we present another way to control the
function $c_{k}$ in a numerical context.\ \ \

\section{Poisson like distribution}

\noindent We still consider the function $c_{k}$ given by:
\[
c_{k}=\sum \limits_{n=1}^{\infty }\frac{\mu ( n) }{n^{\alpha }}{\left(
1-\frac{1}{n^{\beta }}\right) }^{k}
\]

\noindent Then
\begin{align*}
c_{k}&=\sum \limits_{n=1}^{\infty }\frac{\mu ( n) }{n^{\alpha }}
e ^{k \ln ( 1-\frac{1}{n^{\beta }}) }
\\%
 &=\sum \limits_{n=1}^{\infty }\frac{\mu ( n) }{n^{\alpha }} e ^{-\frac{k}{n^{\beta
}}} e ^{k( \ln ( 1-\frac{1}{n^{\beta }}) +\frac{1}{n^{\beta }})
}
\\%
 &=\sum \limits_{n=1}^{\infty }\frac{\mu ( n) }{n^{\alpha }} e ^{-\frac{k}{n^{\beta
}}} e ^{\Delta ( k,n,\beta ) }
\end{align*}

\noindent Notice that $\Delta <0$. For $\beta$ large we set $\Delta
=0$ to obtain the following approximation: 
\[
c_{k}\cong \sum \limits_{n=1}^{\infty }\frac{\mu ( n) }{n^{\alpha
}} e ^{-\frac{k}{n^{\beta }}}
\]

\noindent In this approximation we see that $c_{k}$ of (8) becomes
equal to $\psi ( x=k) $ of (1). Moreover: 
\begin{align*}
c_{k}&\cong \sum \limits_{n=1}^{\infty }\frac{\mu ( n) }{n^{\alpha
}} e ^{k( 1-\frac{1}{n^{\beta }}) } e ^{-k}
\\%
 &=\sum \limits_{n=1}^{\infty }\frac{\mu ( n) }{n^{\alpha }}\sum
\limits_{p=0}^{\infty }\frac{k^{p}}{p!}{\left( 1-\frac{1}{n^{\beta
}}\right) }^{p} e ^{-k}
\end{align*}

\noindent Thus:
\begin{equation}
c_{k}\cong \sum \limits_{p=0}^{\infty }c_{p}\frac{k^{p}}{p!} e ^{-k}
\end{equation}

\noindent We are in the presence of a Poisson distribution for the
$c_{p}$: in this way, in our numerical computations, we may control
in a ``more satisfactory'' way the values of $c_{k}$. The approximation
for $c_{k}$ by means of the Poisson distribution for the $c_{k}$`s
we found, will be more satisfactory with increasing values of $\beta$
and for large values of $k$. We will also use the approximation
given by (9) in which the upper limit of summation will be given
by $2k$ instead of $ \infty $, i.e. for large {\itshape k},
\begin{equation}
c_{k}\cong \sum \limits_{p=0}^{2k}c_{p}\frac{k^{p}}{p!} e ^{-k}
\end{equation}

\section{Numerical experiments}

\subsection{The case $\alpha =\frac{7}{2}$ and $\beta = 4$}

\noindent This is a case of interest since the behaviour of the
$c_{k}$ at large values of $k$ is expected to be the same as the
case $\alpha =\beta =2$ [2, 8]. In fact from (7) we ask that for
$\mathfrak{R}( s) >\frac{1}{2}$: 
\begin{equation}
\left| c_{k}( 7/2,4) \right| \leq \frac{C}{k^{\frac{7/2-1/2-\epsilon
}{4}}}\sim \frac{k^{\frac{\epsilon }{4}}}{k^{\frac{3}{4}}}\sim \left|
c_{k}( 2,2) \right| 
\end{equation}

\noindent In the figures 1, 2 and 3 we give the plot respectively
of $\log |c_{k}|,\log ( |c_{k}\log  k|) $ and $\log ( |{c_{k}( \log
k) }^{2}|) \text{}$ as a function of $\log  k$ for $k$ up to 1000
together with the straight line with slope $-\frac{3}{4}$ which
is tangent to the curves at some point. The $c_{k }$ were computed
calculating (4) until $n=10000$.

\noindent This experiment indicates that $c_{k}$ decays more fast
than $\frac{C}{{(\log  k)}^{2}k^{\frac{3}{4}}}$ as announced by
Baez-Duarte in [3] for the case $\alpha =\beta =2$, i.e. more fast
then the bound (11) if the RH is true (see the necessary condition
in Appendix B).\ \ \ 
\begin{figure}[h]
\begin{center}
\includegraphics{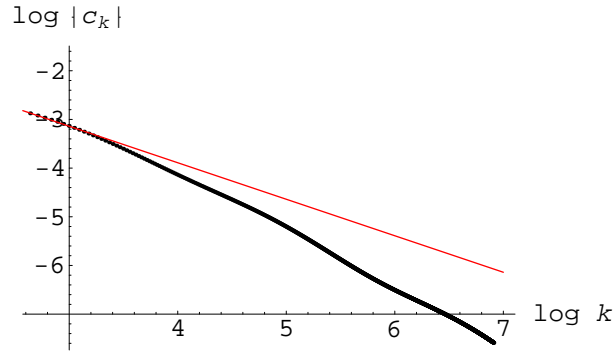}
\caption{Plot of $\log |c_{k}|=C-\frac{3}{4}\log  k$
together with the straight line of slope $-\frac{3}{4}$.}
\end{center}
\end{figure}

\begin{figure}[h]
\begin{center}
\includegraphics{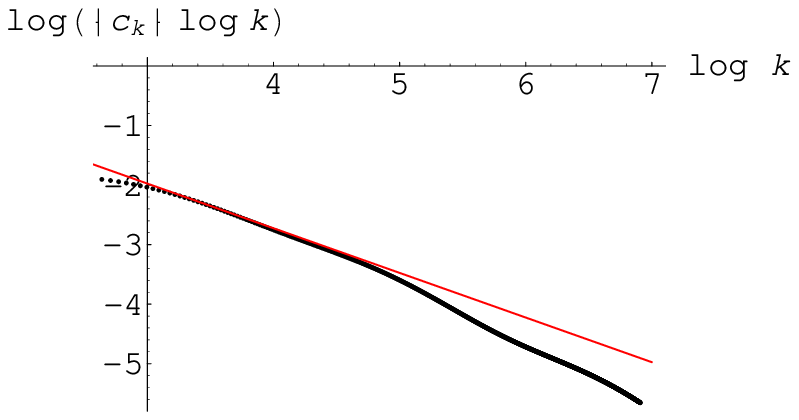}
\caption{Plot of $\log ( |c_{k}|\log  k) =C-\frac{3}{4}\log
k$ together with the tangent straight line of slope $-\frac{3}{4}$.}
\end{center}
\end{figure}

\begin{figure}[h]
\begin{center}
\includegraphics{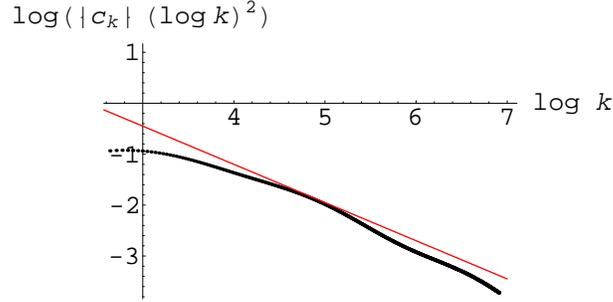}
\caption{Plot of $\log ( |c_{k}|{(\log  k)}^{2}) =C-\frac{3}{4}\log
k$ together with the tangent straight line of slope $-\frac{3}{4}$.}
\end{center}
\end{figure}

\subsection{The case $\alpha =\frac{2+3\beta}{4}$}

\noindent If $\alpha = \frac{2+3\beta
}{4}$ then all the $c_{k}( \alpha ,\beta ) $ are expected to have
a decay similar to $k^{-\frac{3}{4}}$ (7). Below (see figure 4)
we present the numerical results for $\beta$ from 1 to 6. This
is in agreement with the theory of Appendix B. 
\begin{figure}[h]
\begin{center}
\includegraphics{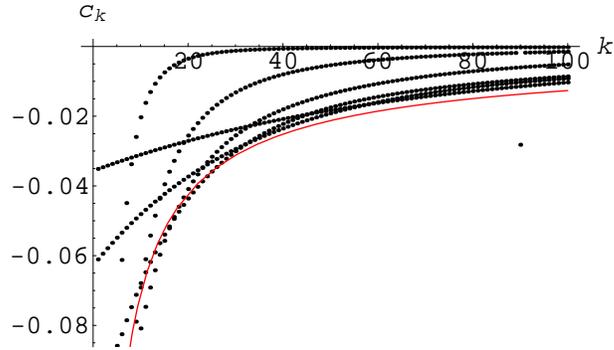}
\caption{Plot of $c_{k}( \frac{2+3\beta }{4},\beta
) $ for $\beta =1, 2,3 ,4, 5, 6$ (dotted lines) vs. the reference
curve $f( k) =-0.4 k^{-\frac{3}{4}}$.}
\end{center}
\end{figure}

\subsection{The case $\alpha =\frac{7}{2}$ with $\beta \rightarrow \infty$}

\noindent Let $\alpha =\frac{7}{2}$ be fixed, from (8) as $\beta$
increases we get:
\[
\begin{array}{rl}
 \operatorname*{\lim }\limits_{\beta \,\rightarrow \:\infty }c_{k}
& =\operatorname*{\lim }\limits_{\beta \,\rightarrow \:\infty }\sum
\limits_{j=0}^{k}{\left( -1\right) }^{j}\binom{k}{j}\frac{1}{\zeta
( 7/2+\beta j) }=\binom{k}{0}\frac{1}{\zeta ( 7/2) }-1+1+\sum \limits_{j=1}^{k}\binom{k}{j}{\left(
-1\right) }^{j} \\
  & =\frac{1}{\zeta ( 7/2) }-1+\sum \limits_{j=0}^{k}\binom{k}{j}{\left(
-1\right) }^{j}
\end{array}
\]

\noindent Thus:
\begin{equation}
\operatorname*{\lim }\limits_{\beta \,\rightarrow \:\infty }c_{k}=\frac{1}{\zeta
( 7/2) }-1\cong - 0.112479\ \ \ \ \ \ \ \ \ \ \ \ \ \ \ \ \ \ \ \ \ \ \ \ \ \forall
k\in \mathbb{N}
\end{equation}
\begin{figure}[h]
\begin{center}
\includegraphics{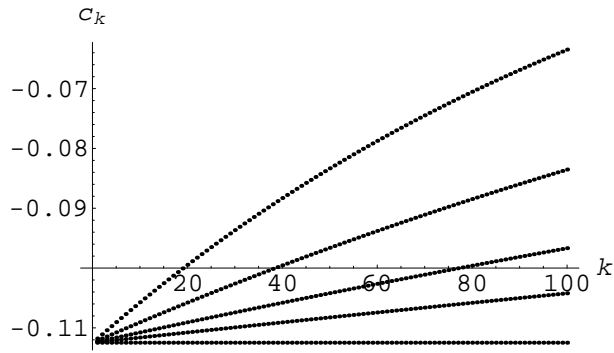}
\caption{Plot of $c_{k}$ for $\alpha =7/2$ and $\beta
=4,5,6,7,20$ (from top to bottom).}
\end{center}
\end{figure}

\noindent Our numerical experiments convalidate these results. We
calculated the first 100 $c_{k}$ for $\beta =4,5,6,7,20$. For $\beta
=20$ we get already a convergence to the theoretical limit (12),
see figure 5. So, this infinite $\beta$ limit obtained by the numerical
calculations for low values of {\itshape k} (up to 100) indicates
that RH may barely be true (see (7) as $\beta \rightarrow \infty
$).

\subsection{The Poisson distribution}

\noindent To demonstrate the goodness of the approximation's formula
(10) we computed the $c_{k}$ until $k=1000$ for the case $\alpha
=\frac{7}{2}, \beta =4$ (using (4)). Then using these already computed
$c_{k}$ we calculated also the first 500 $c_{k}$ of (10). We ploted
these two curves together. In figure 6 we see that from $k\cong
40$ the Poisson approximation is essentially the same as the real
function.\ \ \ \ 
\begin{figure}[h]
\begin{center}
\includegraphics{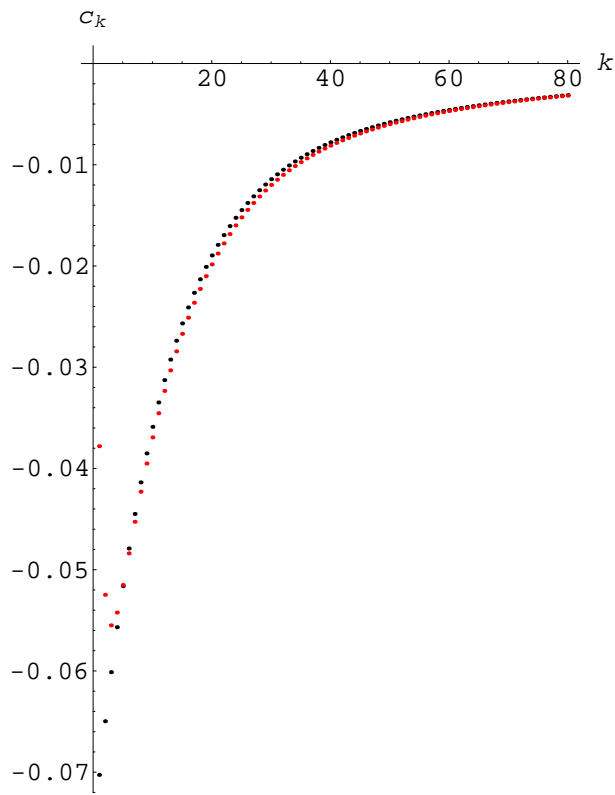}
\caption{Plot of $c_{k}$ for $\alpha =\frac{7}{2},
\beta =4$ (left) vs. the Poisson approximation (right).}
\end{center}
\end{figure}

\subsection{A case with different values of $\rho$ in the critical strip}

\noindent This experiment is carried out in the critical strip in
order to compare the behaviour of $c_{k}$ in three different cases
of $\rho$. In fact from Appendix B the $c_{k}$ should all decay
at least as:
\begin{equation}
\left| c_{k}\right| <\frac{A}{k^{\frac{3}{4}-\epsilon }}
\end{equation}

\noindent In particular for the cases we treat, i.e. $\alpha =\beta
=2$ (the case of Riesz at $\rho =\frac{1}{2}$), $\alpha =\beta =3$
at $\rho =\frac{3}{4}$ and $\alpha =\beta =4$ at $\rho =1$, the
results are presented on figure 7, 8 and 9 respectively. The plots
indicate similar behaviour not in disagreement with (13). Further,
by looking at these plots, the numerical computations visualize
that: 
\[
\left| c_{k}( 2,2) \right| \leq \left| c_{k}( 3,3) \right| \leq
\left| c_{k}( 4,4) \right| 
\]

\noindent We know that $|c_{k}( 4,4) |\ll k^{-\frac{3}{4}+\epsilon
}$ for $\rho \geq 1$ and so $|c_{k}( 3,3) |\ll k^{-\frac{3}{4}+\epsilon
}$. This would indicate that there is no zero for $\rho >\frac{3}{4}$.
Finally $|c_{k}( 2,2) |\ll k^{-\frac{3}{4}+\epsilon }$ would confirm
that there is also no zero for $\rho >\frac{1}{2}$.\ \ \ \ 

\begin{figure}[h]
\begin{center}
\includegraphics{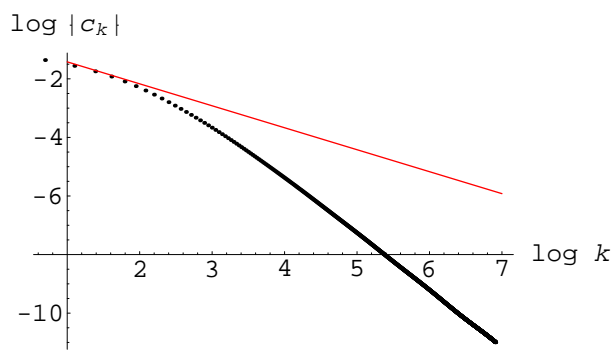}
\caption{Plot of $\log ( |c_{k}|) $ for $\alpha =\beta
=2$ together with the straight line of slope $-\frac{3}{4}$.}
\end{center}
\end{figure}

\begin{figure}[h]
\begin{center}
\includegraphics{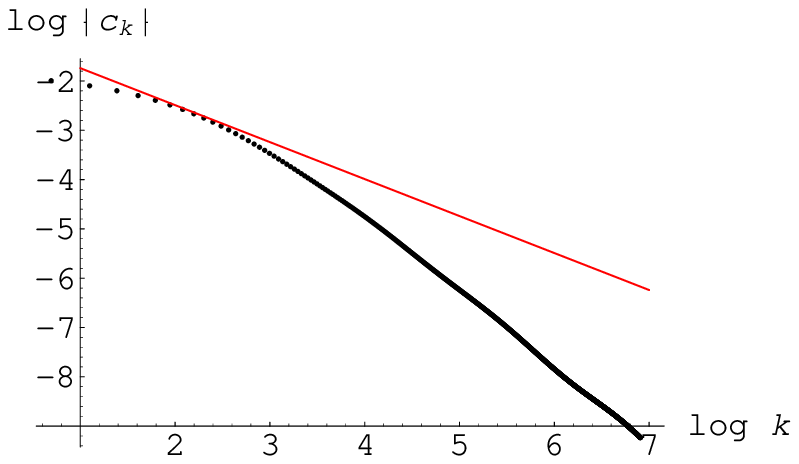}
\caption{Plot of $\log ( |c_{k}|) $ for $\alpha =\beta
=3$ together with the straight line of slope $-\frac{3}{4}$.}
\end{center}
\end{figure}

\begin{figure}[h]
\begin{center}
\includegraphics{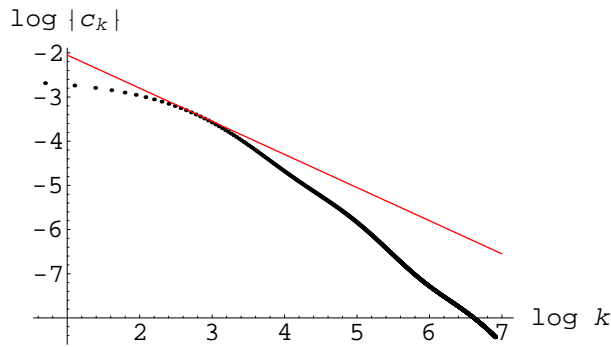}
\caption{Plot of $\log ( |c_{k}|) $ for $\alpha =\beta
=4$ together with the straight line of slope $-\frac{3}{4}$.}
\end{center}
\end{figure}

\subsection{The case $\alpha =\frac{3}{2}$ and $\beta = 1$}

\noindent In this case, from the Riesz formulation (1) or from (6),
${[\zeta ( s) ]}^{-1}$ is different from infinity if $c_{k}$ decays
at least as $k^{-1}$. Setting $x=k$ the Riesz function (1) is given
by: 
\[
f( k) =\sum \limits_{n=1}^{\infty }\frac{\mu ( n) }{n^{3/2}}e^{-\frac{k}{n}}
\]

\noindent Below (figure 10) we present the plot of $f( k) $ together
with that of $g( k) =-0.07 k^{-1}$ for {\itshape k} up to 600 and
for a maximum value of only some hundred for {\itshape n}.
\begin{figure}[h]
\begin{center}
\includegraphics{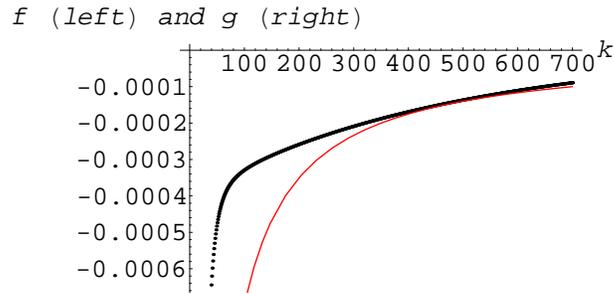}
\caption{Plot of {\itshape f} and {\itshape g} ($\alpha
=\frac{3}{2}, \beta =1$).}
\end{center}
\end{figure}

\subsection{The case $\alpha =\frac{1}{2}$}

\noindent In this case ($\alpha <1$!) we cannot employ the argument
of Appendix A, but we have for $\mathfrak{R}( s) -\epsilon \geq
\rho =\frac{1}{2}$:
\[
\left| \varphi ( s) \right| <\sum \limits_{k=0}^{\infty }k^{-\left(
\frac{\epsilon }{\beta }+1\right) }\left| c_{k}( 1/2,\beta ) \right|
\leq \sum \limits_{k=0}^{\infty }k^{-\left( \frac{\epsilon }{\beta
}+1\right) }\left| c_{k}( 1/2,\infty ) \right| 
\]

\noindent From Subsection 4.3 we know that $|c_{k}( 1/2,\infty )
|=|\frac{1}{\zeta ( 3/2) }-1|\cong 0.617$, thus for any finite $\beta$,
$\varphi ( s) $ is also finite.

\section{Conclusions}

\noindent In this work we have revisited the original criteria of
Riesz and of Hardy and Littlewood for the Riemann Hypothesis in
light of recent pioneering works concerning the possible representations
of the Riemann Zeta function by means of the Pochammer's polynomials.
The discrete representations in the case $\alpha =\beta =2$ are
do to Maslanka and to Baez-Duarte. In order to carry out our numerical
experiments related to the criteria, we have first extended the
analytical formulation to a more general class of functions containing
two parameters $\alpha$ and $\beta$; using a theorem of Baez-Duarte
we have specified a sufficient and necessary condition for the truth
of the RH for our general class of functions i.e. for the decay
of the coefficients $c_{k}$ as a power law of {\itshape k}. Moreover
in doing this we have found the emergence of a Poisson-like distribution
for the $c_{k}$ which should be exact in the large $\beta$ limit.
Numerical experiments have been carried out for various cases.
\begin{enumerate}
\item For $\alpha =\frac{7}{2}$ and $\beta =4$ we have presented
intensive calculation using the M\"obius function up to $n=10000$
and for {\itshape k} up to some hundreds. For this case, the power
law decay $k^{-\frac{3}{4}}$ is the same as that appearing in the
original work of Riesz ($\alpha =\beta =2$) and also investigated
numerically by Baez-Duarte. The experiments confirm the correctness
of the power law within the range of the values of {\itshape n}
and of {\itshape k} we were able to treat.
\item For $\alpha$ and $\beta$ such that the $c_{k}$ should all
give the power law decay $k^{-\frac{3}{4}}$ at large values of {\itshape
k} to ensure the truth of the RH, i.e those where $\alpha =\frac{2+3\beta
}{4}$, we have presented some experiments at low values of {\itshape
k} and some values of $\beta$ which confirms this power law decay.
All functions $c_{k}$ assume negative values with plots lying above
a fixed curve of equation $y=A k^{-\frac{3}{4}}$ for some fixed
constant {\itshape A }independent of $\beta$. In addition for this
expected behaviour we have reported the results of intensive computations
for the cases $\alpha =\beta =2$, $\alpha =\beta =3$ and\ \ $\alpha
=\beta =4$. 
\item Finally we have considered some experiments in the large $\beta$
limit which indicate that the plots of $c_{k}$ become more and more
flat, well approximated by the mean value of the Poisson-type distribution
we have founded. As $\beta$ becomes large and large the $c_{k}$
approaches a constant value, for all {\itshape k}, indicating that
in this sense the RH may barely be true.
\end{enumerate}

\noindent This work, still accompanied by numerical experiments,
may be expanded in the search of other new representations of the
Riemann Zeta function, different of the one considered here; moreover
there is the aim that the new criteria will be useful in the context
of more numerical experiments. The works will be presented in a
near future [4, 5].

\section*{Appendix A}
\noindent We follow strictly the lines of calculations of Baez-Duarte
[2] to show that the representation (6) for ${[\zeta ( s) ]}^{-1}$
is unconditionally valid for $\mathfrak{R}( s) >\rho =1, \alpha
>1$ and $\beta >0$. We consider the quantity:
\[
q_{k}=\sum \limits_{n=1}^{\infty }\frac{1}{n^{\alpha }}{\left( 1-\frac{1}{n^{\beta
}}\right) }^{k}
\]

\noindent Using the MacLaurin series (restricting ourselves to the
main contribution), we have that:
\[
q_{k}\cong \int _{1}^{\infty }\frac{1}{x^{\alpha }}{\left( 1-\frac{1}{x^{\beta
}}\right) }^{k}dx
\]

\noindent Then with the variable change $y=\frac{1}{x^{\beta }}$
we obtain:
\begin{align*}
q_{k}&\cong \frac{1}{\beta }\int _{0}^{1}{y^{\frac{\alpha -1}{\beta
}-1}( 1-y) }^{k+1-1}dy
\\%
 &=B( \frac{\alpha -1}{\beta },k+1) 
\end{align*}

\noindent where
\[
B( \lambda ,\mu ) =\int _{0}^{1}{x^{\lambda -1}( 1-x) }^{\mu -1}dx=\frac{\Gamma
( \lambda ) \Gamma ( \mu ) }{\Gamma ( \lambda +\mu ) }
\]

\noindent is the Beta function.

\noindent Thus for $k$ large, we have
\[
q_{k}\cong \frac{1}{\beta }\Gamma ( \frac{\alpha -1}{\beta }) C\frac{k}{k^{\frac{\alpha
-1}{\beta }+1}}\cong \frac{1}{k^{\frac{\alpha -1}{\beta }}}
\]

\section*{Appendix B}

\noindent Still following Baez-Duarte [2] and here for the family
of functions with parameters $\alpha$ and $\beta$, we now show the
necessity of the condition (7), assuming the RH to be true in the
seminfinite strip $\mathfrak{R}( s) >\rho =\frac{1}{2}$. 

\noindent We set $M( x) =\sum \limits_{n\leq x}^{ }\mu ( n) $, then
we obtain $\forall \epsilon >0$:
\[
M( x) \leq x^{\rho +\epsilon }
\]

\noindent Integration by parts gives for the main contribution:
\[
c_{k}=\int _{1}^{\infty }M( x) \frac{d}{dx}\left( \frac{1}{x^{\alpha
}}{\left( 1-\frac{1}{x^{\beta }}\right) }^{k}\right) dx
\]

\noindent With the variable change $y=\frac{1}{x}$, using $M( \frac{1}{y})
\ll y^{-\rho -\epsilon }$ for $y\downarrow 0$ (RH) we have:
\[
\left| c_{k}\right| <\alpha \int _{0}^{1}{y^{\alpha -\rho -\epsilon
-1}( 1-y^{\beta }) }^{k}dy+\beta k\int _{0}^{1}{y^{\alpha +\beta
-\rho -\epsilon -1}( 1-y^{\beta }) }^{k-1}dy
\]

\noindent and finally with $y^{\beta }=z$ we obtain
\[
\left| c_{k}\right| <\frac{\alpha }{\beta }\int _{0}^{1}{{z^{\frac{\alpha
-\rho -\epsilon }{\beta }-1}( 1-z) }^{k+1-1}}^{ }dz+k\int _{0}^{1}{z^{\frac{\alpha
-\rho -\epsilon +\beta }{\beta }-1}( 1-z) }^{k-1}dz
\]

\noindent which for large $k$ is given by:
\[
\left| c_{k}\right| <\frac{\alpha }{\beta }\frac{\Gamma ( \frac{\alpha
-\rho -\epsilon }{\beta }) }{k^{\frac{\alpha -\rho -\epsilon }{\beta
}}}+\frac{\Gamma ( \frac{\alpha -\rho -\epsilon +\beta }{\beta })
}{k^{\frac{\alpha -\rho -\epsilon +\beta }{\beta }+1}}<\frac{C}{k^{\frac{\alpha
-\rho -\epsilon }{\beta }}}
\]

%\section*{References}
\Bibliography{9}
\item Baez-Duarte L 2003 {\it arXiv:math.NT/0307214v1} 16 July 2003
\item Baez-Duarte L 2003 {\it arXiv:math.NT/0307215v1} 16 July 2003
\item Baez-Duarte L 2005 {\it arXiv:math.NT/0504402v1} 20 April 2005
\item Beltraminelli S and Merlini D 2006 {\it in preparation, to be posted}
\item Beltraminelli S and Merlini D 2006 {\it in preparation, to be posted}
\item D'Errico M \etal 2004 (unpublished) {\it presented at the International Workshop on Complex Systems (Cerfim-Issi)} held in Locarno (Switzerland), 16-18 september 2004
\item Hardy GH and Littlewood JE 1918 {\it Acta Mathematica} {\bf 41} 119
\item Maslanka K 2001 {\it arXiv:math-ph/0105007v1} 4 May 2001
\item Riesz F 1916 {\it Acta Mathematica} {\bf 40} 185

\endbib

\end{document}